\documentclass[12pt]{article}
\usepackage{color}
\usepackage{amssymb, amsmath}
\usepackage{epsfig}
\usepackage{graphicx}

\newtheorem{theorem}{Theorem}

\newtheorem{lemma}{Lemma}%

\newcommand{\proof}{\noindent\textbf{Proof.~}}

\newcommand{\qed}{\space\hfill\hspace*{\fill} $\vbox{\hrule\hbox{\vrule
height1.3ex\hskip1.3ex\vrule}\hrule}$\hss\vskip\topsep\relax}

\begin{document}

\title{An explicit and positivity preserving numerical scheme for the mean reverting CEV model}

\author{Nikolaos Halidias \\
{\small\textsl{Department of Mathematics }}\\
{\small\textsl{University of the Aegean }}\\
{\small\textsl{Karlovassi  83200  Samos, Greece} }\\
{\small\textsl{email: nikoshalidias@hotmail.com}}}

\maketitle

\begin{abstract} In this paper we propose an explicit and positivity preserving scheme
for the mean reverting  CEV model which converges in the mean
square sense with convergence order $a(a-1/2)$.
\end{abstract}



\section{Introduction}
Let $(\Omega, {\cal F}, \mathbb{P}, {\cal F}_t)$ be a complete
probability space with a filtration and let  a Wiener process
$(W_t)_{t \geq 0}$ defined on this space. We consider here the
mean reverting CEV process,
\begin{eqnarray}
x_t = x_0 + \int_0^t (kl - k x_s)ds + \sigma \int_0^t x_s^{a}
dW_s,
\end{eqnarray}
where $k,l,\sigma \geq 0$ and $a \in (1/2,1)$. It is well known
that this sde has a unique strong solution which is strictly
positive. Our starting point was the paper of \cite{Alfonsi3} in
which the author proposes an implicit and positivity preserving
numerical scheme to approximate the above process. This stochastic
process plays important role in financial mathematics. If one
wants to use the above model to price complicate path-dependent
options  maybe it is useful to approximate it numerically. The
usual Euler scheme (see \cite{Gyongy}) does not preserve
positivity and therefore numerical schemes with this feature are
needed. In this direction on can see \cite{Alfonsi3} and
\cite{Berkaoui2}. Our goal here is to present an explicit,
positivity preserving numerical scheme that converges in the mean
square sense to the true solution with, at least, $a(a-1/2)$ order
of convergence. There must be extensive numerical experiments to
compare all these methods and decide which of them is the best in
each set of parameters.

We will construct our scheme  using the semi discrete method that
we have proposed in \cite{Halidias1} (and further extended in
\cite{Halidias2}, \cite{Halidias3}).
 Let $0 = t_0 <
t_1 < ...<t_n = T$ and set $\Delta = \frac{T}{n}$. Consider the
following stochastic process
\begin{eqnarray}
y_t = \left| \sigma(1-a) (W_t - W_{t_k}) +
\left(y_{t_k}(1-k\Delta)+ \Delta (kl - \frac{a\sigma^2
y_{t_k}^{2a-1}}{2}) \right)^{1-a} \right|^{\frac{1}{1-a}} =
|z_t|^{\frac{1}{1-a}},
\end{eqnarray}
for $t \in (t_k,t_{k+1}]$  where
\begin{eqnarray*}
z_t = \sigma(1-a)  (W_t - W_{t_k}) + \left(y_{t_k}(1-k\Delta) +
\Delta (kl - \frac{a\sigma^2 y_{t_k}^{2a-1}}{2}) \right)^{1-a},
\end{eqnarray*}
for $t \in (t_k,t_{k+1}]$. To show that this stochastic process is
well defined we will check whether
\begin{eqnarray*}
y_{t_k}(1-k\Delta)+ \Delta (kl - \frac{a\sigma^2
y_{t_k}^{2a-1}}{2}) \geq 0.
\end{eqnarray*}
For $kl \geq \frac{a \sigma^2}{2}$, $\Delta \leq \frac{1}{k +
\frac{a \sigma^2}{2}}$ and noting that $0 < 2a - 1  < 1 $ we have
\begin{eqnarray*}
\left(y_{t_k}(1-k\Delta)+ \Delta (kl - \frac{a\sigma^2
y_{t_k}^{2a-1}}{2})\right) \mathbb{I}_{ \{ y_{t_k} > 1 \} } +
\left(y_{t_k}(1-k\Delta)+ \Delta (kl - \frac{a\sigma^2
y_{t_k}^{2a-1}}{2})\right) \mathbb{I}_{ \{ y_{t_k} \leq  1 \} } \geq \\
\mathbb{I}_{ \{ y_{t_k} > 1 \} } y_{t_k}^{2a-1} \left(1-k\Delta  -
\frac{\Delta a \sigma^2}{2}\right)  + \Delta \left(kl - \frac{a
\sigma^2}{2}\right) \geq 0
\end{eqnarray*}
Therefore we impose the following assumption.

{\bf Assumption A} We assume that  $x_0 \in \mathbb{R}_+$.
Moreover, we suppose that $$kl \geq \frac{a \sigma^2}{2}, \quad
\Delta \leq \frac{2}{2k + a \sigma^2}.$$ This stochastic process,
using Ito's formula, has the differential form,
\begin{eqnarray*}
y_t = y_{t_k}(1-k\Delta)+ \Delta (kl - \frac{a\sigma^2
y_{t_k}^{2a-1}}{2}) + \int_{t_k}^t \frac{a \sigma^2}{2}
y_{s}^{2a-1}  ds + \sigma \int_{t_k}^t y_s^{a} sgn(z_s) dW_s,
\quad t \in (t_k,t_{k+1}],
\end{eqnarray*}
Concluding, the numerical scheme that we propose for the mean
reverting CEV model is the following,
\begin{eqnarray*}
y_{t_{k+1}} = \left|\sigma(1-a) (W_{t_{k+1}} - W_{t_k}) +
\left(y_{t_k}(1-k\Delta)+ \Delta (kl - \frac{a\sigma^2
y_{t_k}^{2a-1}}{2}) \right)^{1-a} \right|^{\frac{1}{1-a}}.
\end{eqnarray*}

\section{Main results}
\begin{lemma}
Under Assumption A we have the  estimate for the following
probability,
\begin{eqnarray*}
\mathbb{P}(z_t \leq 0) \leq C
\frac{\Delta^{a-\frac{1}{2}}}{e^{\frac{C}{\Delta^{2a-1}}}}.
\end{eqnarray*}
Therefore, the above probability tend to zero faster than any
power of $\Delta$.
\end{lemma}

\proof Indeed, we have
\begin{eqnarray*}
\mathbb{P}(z_t \leq 0) & = & \mathbb{P} \left( W_t - W_{t_{k+1}}
\leq -\frac{ \left(y_{t_k} (1-k\Delta) + \Delta d
\right)^{1-a}}{\sigma(1-a)}  \right) \\ & \leq & \frac{1}{\sqrt{2
\pi}}\int_{C \Delta^{1/2-a}}^{\infty} e^{-\frac{y^2}{2}} dy \\ &
\leq & C \frac{\Delta^{a-1/2}}{e^{\frac{C}{\Delta^{2a-1}}}},
\end{eqnarray*}
where $d = kl - \frac{a \sigma^2}{2}$. To obtain the last
inequality one can use  problem 9.22, p.112 of \cite{karatzas}.
 \qed

Next we will use the following compact form for our scheme,
\begin{eqnarray*}
y_t = x_0 & & + \int_0^t (kl - k y_{\hat{s}}) ds + \int_0^t
\frac{a \sigma^2}{2} ( y_s^{2a-1} - y_{\hat{s}}^{2a-1}) ds
\\ & &- \int_{t}^{t_{k+1}} \left(kl - k y_{t_k} - \frac{a
\sigma^2}{2} y_{t_k}^{2a-1} \right) ds + \sigma \int_0^t
y_s^{a}sgn(z_s)dW_s, \quad t \in (t_k, t_{k+1}].
\end{eqnarray*}

Consider the following process,
\begin{eqnarray*}
v_t = x_0 + Tkl + \int_0^t \frac{a \sigma^2}{2} y_s^{2a-1} ds +
\int_0^t \sigma y_s^{a}sgn(z_s) dW_s,
\end{eqnarray*}
Then it is clear that $0 \leq y_t \leq v_t$. We will show that
$v_t$ has bounded moments and therefore $y_t$ has also bounded
moments.

\begin{lemma}[Moment bounds]
Under Assumption A  we have the moment bounds,
\begin{eqnarray*}
\mathbb{E}y_t^2 + \mathbb{E} x_t^2  < C,
\end{eqnarray*}
for some $C > 0$
\end{lemma}

\proof Consider the stopping time $\theta_R = \inf \{ t \geq 0 :
v_t > R \}$. Using Ito's formula on $v_{t \wedge \theta_R}^2$ we
obtain,
\begin{eqnarray*}
v_{t \wedge \theta_R}^2 = (x_0+Tkl)^2 + \int_0^t \left(
\frac{2a\sigma^2}{2} v_{s \wedge \theta_R} y_{s \wedge
\theta_R}^{2a-1} + \frac{\sigma^2}{2} y_{s \wedge \theta_R}^{2a}
\right) ds + 2 \sigma \int_0^t v_{ s \wedge \theta_R} y_{s \wedge
\theta_R}^{a}sgn(z_{s \wedge \theta_R})dW_s.
\end{eqnarray*}
Taking expectations on both  sides and noting that $y_t \leq v_t$,
we arrive at
\begin{eqnarray*}
\mathbb{E}v_{t \wedge \theta_R}^2 & \leq & \mathbb{E} (x_0+Tkl)^2
+ C \int_0^t \left( \mathbb{E} v_{s \wedge \theta_R}^2 \right)^{a}
ds
\end{eqnarray*}
Using now a Gronwall type theorem  (see \cite{Mitrinovic}, Theorem
1, p. 360), we arrive at
\begin{eqnarray}
\mathbb{E} v_{t \wedge \theta_R}^2 \leq
 C
\end{eqnarray}
 But $\mathbb{E} v_{t \wedge
\theta_R}^2 = \mathbb{E} (v_{t \wedge \theta_R}^2 \mathbb{I}_{ \{
\theta_R \geq t \} }) +
 R^p P ( \theta_R < t  )$. That means that $P( t \wedge \theta_R < t) = P ( \theta_R < t  ) \to 0$ as
 $R \to \infty$ so $t \wedge \theta_R \to t$ in probability and noting that $\theta_R$ increases as $R$ increases we have that $t \wedge \theta_R \to t$
  almost surely too, as $R \to \infty$. Going back to (3) and
 using  Fatou's lemma  we obtain,
 \begin{eqnarray*}
\mathbb{E} v_t^2 \leq C
\end{eqnarray*}
 The same holds for $x_t$.
 \qed

Next we define the process
\begin{eqnarray*}
h_t = x_0 + \int_0^t (kl - k y_{\hat{s}}) ds + \sigma \int_0^t
y_s^{a}sgn(z_s)dW_s
\end{eqnarray*}
We will show that $h_t,y_t$ remain close.

\begin{lemma}
We have the following estimates,
\begin{eqnarray*}
\mathbb{E}|y_s - y_{\hat{s}}|^2 & \leq & C \Delta \mbox{ for any } s \in [0,T] \\
\mathbb{E}|h_s-y_s|^2 & \leq & C \Delta^{2a-1} \mbox{ for any } s
\in [0,T]
 \\
\mathbb{E} |h_s - y_{\hat{s}}|^2 & \leq & C \Delta^{2a-1}   \mbox{
for any  } s
\in [0,T]   \\
\mathbb{E}|h_s|^2 & < &  A,   \mbox{ for any } s \in [0,T].
\end{eqnarray*}
\end{lemma}

\proof Using the moment bound for $y_t$ we  easily obtain the fact
that
\begin{eqnarray*}
\mathbb{E}|y_s - y_{\hat{s}}|^2 \leq C \Delta
\end{eqnarray*}
and then
\begin{eqnarray*}
\mathbb{E}|h_s - y_s |^2 \leq C \Delta^{2a-1}.
\end{eqnarray*}

Next, we have
\begin{eqnarray*}
\mathbb{E}|h_s - y_{t_k}|^2 \leq 2 \mathbb{E} |h_s - y_s|^2 + 2
\mathbb{E}|y_s-y_{t_k}|^2 \leq  C \Delta^{2a-1}.
\end{eqnarray*}

Finally, to get the moment bound for $h_t$ we just use the fact
that is close to $y_t$, i.e.
\begin{eqnarray*}
\mathbb{E}h_t^2 \leq 2\mathbb{E} |h_t-y_t|^2 +  2 \mathbb{E} y_t^2
\leq C.
\end{eqnarray*}

 \qed

\begin{lemma}[Inverse Exponential Moments]
For the true solution $x_t$ it holds the following bound,
\begin{eqnarray*}
\mathbb{E} \left(\exp(\frac{C}{x_t^{2(1-a)}}) \right) < \infty,
\end{eqnarray*}
for any $C > 0$.
\end{lemma}

\proof We first transform our equation with $z_t = x_t^{2(1-a)}$.
Then using Ito's formula we deduce that
\begin{eqnarray*}
z_t = z_0 + \int_0^t \left(\sigma^2 (1-a)(1-2a) - 2(1-a)k z_s +
\frac{2(1-a)}{z_s^{\frac{2a-1}{2-2a}}} \right) ds + 2(1-a) \sigma
\int_0^t \sqrt{z_s} dW_s,
\end{eqnarray*}
and if we denote by $b(z) = \sigma^2 (1-a)(1-2a) - 2(1-a)k z +
\frac{2(1-a)}{z^{\frac{2a-1}{2-2a}}}$ then it is easy to see that
for any $M > 0$ there exists a $c_M$ such that $b(z) \geq M - c_m
z$. We construct now the following CIR process
\begin{eqnarray*}
f_t = x_0 + \int_0^t (M- c_M f_s)ds + 2(1-a)\sigma \int_0^t
\sqrt{f_s} dW_s.
\end{eqnarray*}
Using a comparison theorem for stochastic differential equations
(see \cite{karatzas}, prop. 5.2.18) we deduce that $z_t \geq f_t
> 0$ (choosing big enough $M > 0$ in order $f_t$ to be strictly
positive and also the inverse moment bound holds for $f_t$) and
using the inverse exponential moments of \cite{Kuznetsov} we have
the desired result. \qed

\begin{theorem}
Under Assumption A we have
\begin{eqnarray*}
\mathbb{E} | x_t - y_t |^2 \leq C \Delta^{2a(a-1/2)},
\end{eqnarray*}
and therefore the order of convergence is at least $a(a-1/2)$.
\end{theorem}

\proof Using Ito's formula on $|x_{\rho} - y_{\rho}|^2$ for some
stopping time $\rho$, we obtain
\begin{eqnarray*}
\mathbb{E}|x_{\rho} - h_{\rho}|^2 & = & \int_0^{\rho} 2k
\mathbb{E}(x_s - h_s)(x_s-y_{\hat{s}}) + \sigma^2
\mathbb{E}(x_s^{a} - y_s^{a} sgn(z_s))^2 ds
\\ & \leq & \int_0^{\rho} (2k \mathbb{E}|x_s - h_s|^2 + 2k
\mathbb{E}|x_s-h_s||h_s - y_{\hat{s}}| +\sigma^2
\mathbb{E}(x_s^{a} - y_s^{a} sgn(z_s))^2 ds.
\end{eqnarray*}
But
\begin{eqnarray*}
\mathbb{E} (x_s^{a} - y_s^{a} sgn(z_s))^2 & \leq & 2
\mathbb{E}(x_s^{a} - y_s^{a})^2 + 2 \mathbb{E}y_s^{2a}
(1-sgn(z_s))^2 \\ &=& 2\mathbb{E}(x_s^{a} - y_s^{a})^2 + 8
\mathbb{P}(z_t < 0) \mathbb{E} \left(y_s^{2a} | (z_t < 0) \right)
\\ & \leq & 2\mathbb{E}(x_s^{a} - y_s^{a})^2 + C \frac{\Delta^{a-1/2}}{e^{\frac{C}{\Delta^{2a-1}}}} \\ & \leq & C\Delta^{2a(a-1/2)}+ 2\mathbb{E}(x_s^{a} -
h_s^{a})^2.
\end{eqnarray*}
We have used the fact that
\begin{eqnarray*}
\mathbb{E}|y_s^{a}- h_s^{a}|^2 \leq C\mathbb{E}|y_s-h_s|^{2a} \leq
C \left(\mathbb{E} |y_s-h_s|^2 \right)^{a} \leq C
\Delta^{2a(a-1/2)}.
\end{eqnarray*}
Moreover, by Young inequality,
\begin{eqnarray*}
\mathbb{E}|x_s-h_s||h_s - y_{\hat{s}}| \leq
1/2\mathbb{E}|x_s-h_s|^2 + 1/2\mathbb{E}|h_s-y_{\hat{s}}|^2 \leq C
\Delta^{2a-1} + 1/2 \mathbb{E}|x_s-h_s|^2.
\end{eqnarray*}
Therefore,
\begin{eqnarray*}
\mathbb{E}|x_{\rho} - h_{\rho}|^2 \leq C \Delta^{2a(a-1/2)} +
\int_0^{\rho}\left(3k\mathbb{E}|x_s-h_s|^2+ 2\sigma^2
\mathbb{E}(x_s^{a} - h_s^{a})^2 \right)ds
\end{eqnarray*}

Setting now $\gamma_t = \int_0^t \frac{8a^2\sigma^2}{x_s^{2(1-a)}+
y_s^{2(1-a)}} ds$ and using the inequality
\begin{eqnarray*}
|x^{a}-y^{a}||x^{1-a}+y^{1-a}| \leq 2a | x- y|
\end{eqnarray*}
we have
\begin{eqnarray}
\mathbb{E}|x_{\rho} - h_{\rho}|^2 \leq C \Delta^{2a(a-1/2)} +
\int_0^{\rho}  \mathbb{E}|x_s - y_s|^2(3ks+ \gamma_s)^{'}ds
\end{eqnarray}
Define the stopping time
\begin{eqnarray*}
\tau_l = \inf \{ s \in [0,T]: 3ks+ \gamma_s \geq l \}.
\end{eqnarray*}
Now, for $\rho = \tau_l$, we use the change of variables setting
$u = 3ks+ \gamma_s$ and therefore $s = \tau_u$ obtaining
\begin{eqnarray*}
\mathbb{E} (x_{\tau_l}-h_{\tau_l})^2 \leq C \Delta^{2a(a-1/2)} +
\int_0^l \mathbb{E} |x_{\tau_u}-h_{\tau_u}|^2 du.
\end{eqnarray*}
Using Gronwall's inequality we obtain,
\begin{eqnarray}
\mathbb{E}|x_{\tau_l} - h_{\tau_l}|^2 \leq C e^{l}
\Delta^{2a(a-1/2)}.
\end{eqnarray}
Going back to (4), for $\rho = t \in [0,T]$, we have under the
change of variable $u =3ks+ \gamma_s$,
\begin{eqnarray}
\mathbb{E}(x_t - h_t)^2 & \leq & C\Delta^{2a(a-1/2)} +\mathbb{E}
\int_0^{ 3kT+\gamma_T} |x_{\tau_u}-h_{\tau_u}|^2 du \nonumber \\
& \leq & C \Delta^{2a(a-1/2)} + \int_0^{\infty} \mathbb{E} \left(
\mathbb{I}_{ \{ \gamma_T \geq u \} } |x_{\tau_u} - h_{\tau_u}|^2
\right) du.
\end{eqnarray}
Using (5) and noting that
\begin{eqnarray*}
\mathbb{E} \left( \mathbb{I}_{ \{ \gamma_T \geq u \} } |x_{\tau_u}
- h_{\tau_u}|^2 \right) = \mathbb{P}(\gamma_T \geq u) \mathbb{E}
\left( |x_{\tau_u} - h_{\tau_u}|^2 | (\gamma_T \geq u) \right)
\end{eqnarray*}
we arrive at
\begin{eqnarray*}
\mathbb{E}(x_t - h_t)^2 \leq C\Delta^{2a(a-1/2)}\left(1+
\int_0^{\infty} \mathbb{P} (  \gamma_T \geq u) e^{u} du \right).
\end{eqnarray*}
Note that,
\begin{eqnarray*}
\mathbb{P}( \gamma_T \geq u) \leq \frac{1}{e^{u}}
\mathbb{E}(e^{\gamma_T} ),
\end{eqnarray*}
therefore, using the inverse exponential moments of the true
solution we deduce that
\begin{eqnarray*}
\mathbb{E}(x_t - h_t)^2 \leq C\Delta^{2a(a-1/2)}
\end{eqnarray*}
 \qed


\begin{thebibliography}{99}




\bibitem{Alfonsi3} {A. Alfonsi}, {\em Strong order one convergence of a drift implicit Euler
scheme: Application to the CIR process}, Statistics and
Probability Letters, Volume 83, Issue 2, (2013), pp. 602-607.

\bibitem{Andersen} {L. Andersen}, {\em Simple and efficient simulation of the Heston stochastic volatility
model}, Journal of Computational Finance, (2008), Vol. 11, No. 3.


\bibitem{Berkaoui2} {A. Berkaoui, M. Bossy and  A. Diop}, {\em Euler
scheme for SDEs with non-Lipchitz diffusion coefficient: strong
convergence}, ESAIM 12, (2008), 1-11.




\bibitem{Gyongy} {I. Gyongy and M. Rasonyi}, {\em A note on Euler approximations for SDEs with Holder
continuous diffusion coefficients}, Stochastic Processes and their
Applications 121 (2011) 2189-2200.


\bibitem{Halidias1} {N. Halidias}, {\em Semi-discrete   approximations for stochastic differential
equations and applications}, International Journal of Computer
Mathematics, (2012), pp. 780-794.


\bibitem{Halidias2} {N. Halidias}, {\em Construction of positivity preserving numerical schemes for a class of  multidimensional stochastic differential
equations}, Discrete and Continuous Dynamical Systems, 2014.


\bibitem{Halidias3} {N. Halidias}, {\em A novel approach to construct numerical methods for stochastic differential
equations}, Numerical Algorithms May 2014, Volume 66, Issue 1, pp
79-87.




\bibitem{karatzas}{I. Karatzas and S. Shreve}, {\em Brownian Motion and
Stochastic Calculus}, Springer, 1991.


\bibitem{Kuznetsov} {T. R. Hurd and A. Kuznetsov},  {\em Explicit formulas for Laplace transforms of stochastic integrals}, Markov Process. Relat. Fields, 14,
277-290 (2008).




\bibitem{Mitrinovic}{D. S. Mitrinovic, J. E. Pecaric and A. M. Fink},
{\em Inequalities Involving Functions and Their Integrals and
Derivatives}, Kluwer, 1991.




\bibitem{Yamada} { T. Yamada and S.  Watanabe} {\em On the uniqueness of
solutions of stochastic differential equations}. { J. Math. Kyoto
Univ.} \textbf{11,} 155-167, 1971.



\end{thebibliography}
\end{document}